\documentclass[a4paper,12pt]{amsart}
\usepackage{amsfonts,amsmath,amsthm,amssymb}
\usepackage{latexsym}
\usepackage{graphics}
\newtheorem{theorem}{Theorem}[section]

\newtheorem{proposition}[theorem]{Proposition}

\theoremstyle{definition}

\theoremstyle{remark}
\newtheorem{rem}[theorem]{Remark}
\theoremstyle{remark}

\newcommand{\beql}[1]{\begin{equation}\label{#1}}
\newcommand{\eeq}{\end{equation}}

\begin{document}

\title{On quaternary complex Hadamard matrices of small orders}

\author{Ferenc Sz\"oll\H{o}si}


\date{May, 2010., Preprint.}

\address{Department of Mathematics and its Applications, Central European University, H-1051, N\'ador u. 9, Budapest, Hungary}\email{szoferi@gmail.com}

\thanks{This work was supported by the Hungarian National Research Fund OTKA-K77748.}

\begin{abstract}
One of the main goals of design theory is to classify, characterize and count various combinatorial objects with some prescribed properties. In most cases, however, one quickly encounters a combinatorial explosion and even if the complete enumeration of the objects is possible, there is no apparent way how to study them in details, store them efficiently, or generate a particular one rapidly. In this paper we propose a novel method to deal with these difficulties, and illustrate it by presenting the classification of quaternary complex Hadamard matrices up to order $8$. The obtained matrices are members of only a handful of parametric families, and each inequivalent matrix, up to transposition, can be identified through its fingerprint.
\end{abstract}

\maketitle


\section{Introduction}
A complex Hadamard matrix $H$ of order $n$ is an $n\times n$ matrix with unimodular entries satisfying $HH^\ast=nI$ where $\ast$ is the conjugate transpose and $I$ is the identity matrix of order $n$. In other words, any two distinct rows (or columns) of $H$ are complex orthogonal. Complex Hadamard matrices have important applications in quantum optics, high-energy physics \cite{AMM}, operator theory \cite{popa} and in harmonic analysis \cite{MMM}, \cite{tao}. They also play a crucial r\^ole in quantum information theory, for construction of teleportation and dense coding schemes \cite{wer}, and they are strongly related to mutually unbiased bases (MUBs) \cite{karol2}. In finite mathematics and in particular, in design and coding theory, the real Hadamard matrices and their discrete generalization, the $n\times n$ Butson-type matrices \cite{Bu}, composed of $q$-th roots of unity appear frequently. These matrices are the main algebraic objects behind the Fourier and Walsh--Hadamard Transforms and perfect quaternary sequences. They can also be utilized to obtain high performance quaternary error-correcting codes and low correlation $\mathbb{Z}_4$-sequence families \cite{horadam}, and shall be denoted through this paper by $BH(q,n)$. The notations $BH(2,n)$ and $BH(4,n)$ correspond to the real and \emph{quaternary} Hadamard matrices, respectively.

One of the main features of these $BH(q,n)$ matrices is the large symmetry they have, which allows one to introduce some free parameters into them obtaining a continuous family of complex Hadamard matrices. This idea was used in a series of papers \cite{dita}, \cite{dita2}, \cite{dita3}, \cite{MRS}, \cite{SZF}, where new, previously unknown complex Hadamard matrices were constructed in this way. The parametrization of $BH(q,n)$ matrices is in some sense the continuous analogue of the switching operation, a well-known technique in design theory, and allows one to escape equivalence classes and therefore collect many inequivalent matrices into a single parametric family.

The aim of this paper is to illustrate how powerful these continuous families are by presenting the complete characterization of the $BH(4,n)$ matrices up to $n=8$, and describing them in a compact and elegant form in terms of parametric families of complex Hadamard matrices. Interestingly, up to transposition, all equivalence classes can be identified through a recently introduced invariant, the fingerprint \cite{finger}. Constructing these matrices and identifying a member from each equivalence class took only a few days of computer search, and we do not consider it a deep mathematical achievement. The main result of this paper, however, is \emph{the way} we present the representatives of the equivalence classes, as with the aid of only a handful of parametric families of complex Hadamard matrices we can encode a lot of inequivalent $BH(4,n)$ matrices \emph{simultaneously}. In principle, by introducing $m$ free affine parameters into a $BH(q,n)$ matrix one can obtain as many as $q^m$ inequivalent matrices of the same type, therefore describing an \emph{exponential number of} equivalence classes through a single matrix. We believe that the concept of parametrization shall turn out to be a useful tool in discrete mathematics as well, and by utilizing this method one shall be able to obtain new lower bounds on the number of various combinatorial objects in the future.

In the following section we give a full characterization of quaternary complex Hadamard matrices of orders $1$-$8$, and describe them as members of parametric families of complex Hadamard matrices. Throughout the paper we adopt the notation from \cite{karol} for well-known families of Hadamard matrices, such as $F_8^{(5)}(a,b,c,d,e)$, \emph{etc}.

\section{List of quaternary complex Hadamard matrices of small orders}
It is easily seen that quaternary complex Hadamard matrices of order $n$ can exist only if $n=1$ or $n$ is even, and it is conjectured that this condition is sufficient \cite{horadam}. It is trivial to check that there is a unique complex Hadamard matrix of orders $n=1,2$, however, for composite orders starting from $n=4$ one encounters infinite families already. It has been shown that all complex Hadamard matrices of order $4$ are members of the following $1$-parameter family, depending on the unimodular complex number $a$ given in \cite{craigen} as
\[H_4(a)=\left[\begin{array}{cc|cc}
1 & 1 & 1 & 1\\
1 & -1 & a & -a\\
\hline
1 & 1 & -1 & -1\\
1 & -1 & -a & a\\
\end{array}\right].\]
Note that all $2\times 2$ blocks of $H_4(a)$ are complex Hadamard matrices. Matrices arising in this fashion are called \emph{Di\c{t}\u{a}-type} (cf.\ \cite{MRS}). Recall that Hadamard matrices $H$ and $K$ are \emph{equivalent}, if there exists permutation matrices $P_1, P_2$ and unitary diagonal matrices $D_1, D_2$ such that $P_1D_1HD_2P_2=K$ holds. We denote this property by $H\sim K$. Clearly, one can consider matrices containing a full row and column of numbers $1$ due to equivalence. Such matrices are called \emph{dephased}. Let us remark here that the transpose of a Hadamard matrix $H^T$ is also Hadamard, but not necessarily equivalent to $H$.

It is not apparent at first, how many equivalence classes of $BH(4,4)$ matrices are contained in the family $H_4(a)$ for various values of $a$. To quickly check that matrices $H_4(1)$ and $H_4(\mathbf{i})$ are inequivalent, we can utilize Haagerup's invariant $\Lambda(H)$, defined for complex Hadamard matrices $H$ of order $n$, which is the following set \cite{haagerup}:
\[\Lambda(H):=\left\{h_{ij}h_{kl}\overline{h}_{il}\overline{h}_{kj} : i,j,k,l=1,2,\hdots, n\right\}.\]
We readily see that $\Lambda(H_4(1))$ and $\Lambda(H_4(\mathbf{i}))$ are different, therefore these matrices are inequivalent. Also, by exchanging rows $2,4$ one can see that $H(a)\sim H(-a)$ and hence $H(1)\sim H(-1)$ and $H(\mathbf{i})\sim H(-\mathbf{i})$. Clearly the Haagerup invariant corresponding to $BH(q,n)$ matrices can contain the $q$-th roots of unity only.
\begin{rem}\label{myRem}
It might happen that a complex Hadamard matrix $H$ does not look like a $BH(q,n)$ matrix, but is equivalent to one; however, $BH(q,n)$ matrices can be quickly recognized. Indeed, after dephasing $H$ all of its entries should be some $q$-th roots of unity, otherwise $\Lambda(H)$ would contain an entry which is not a $q$-th root of unity showing that $H$ \emph{cannot be} equivalent to a $BH(q,n)$ matrix at all. This implies that checking the number of inequivalent $BH(q,n)$ matrices in dephased parametric families is a \emph{finite process}.
\end{rem}
The above remark and the paragraph just before it easily allow to show that there are precisely $2$ inequivalent $BH(4,4)$ matrices: $H_4(1)$ and $H_4(\mathbf{i})$, respectively.

An interesting $1$-parameter family of complex Hadamard matrices of order $6$ was discovered by Zauner \cite{Z} (and independently, but slightly later by Di\c{t}\u{a} \cite{dita}) who used the $2$-circulant form of the following matrix
\[D_6(c)=\left[
\begin{array}{rrrrrr}
 1 & 1 & 1 & 1 & 1 & 1 \\
 1 & -1 & \mathbf{i} & -\mathbf{i} c & -\mathbf{i} & \mathbf{i} c \\
 1 & \mathbf{i} & -1 & \mathbf{i} c & -\mathbf{i} & -\mathbf{i} c \\
 1 & -\mathbf{i}\overline{c} & \mathbf{i}\overline{c} & -1 & \mathbf{i} & -\mathbf{i} \\
 1 & -\mathbf{i} & -\mathbf{i} & \mathbf{i} & -1 & \mathbf{i} \\
 1 & \mathbf{i}\overline{c} & -\mathbf{i}\overline{c} & -\mathbf{i} & \mathbf{i} & -1
\end{array}
\right]\]
as a seed matrix to exhibit an infinite family of MUB triplets in $\mathbb{C}^6$. A quick computer search amongst all possible $BH(4,6)$ matrices confirms that despite the $1$-parameter freedom we have here, there is only a single equivalence class of $BH(4,6)$ matrices. One might wonder why the matrices $D_6(1)$ and $D_6(\mathbf{i})$ are equivalent. This is not trivial to realize, as one need to enphase one of the matrices to obtain the other one. This particular family of complex Hadamard matrices has been utilized in harmonic analysis as well. Considering the scaled logarithm of the $BH(8,6)$ matrix $D_6((1+\mathbf{i})/\sqrt2)$, one finds a $6\times 6$ integer matrix having rank $3$ over $\mathbb{Z}_8$ leading to the disproval of Fuglede's conjecture in dimension $3$ \cite{MMM}.

Now we turn to the classification of $BH(4,8)$ matrices. An exhaustive computer search revealed after a few days the following

\begin{proposition}\label{myCount}
There are precisely $15$ $BH(4,8)$ matrices, up to equivalence; $5$ of these matrices are equivalent to a symmetric one, while another $10$ are inequivalent to their transpose.
\end{proposition}

Now, having \emph{all} inequivalent $BH(4,8)$ matrices at our disposal, our goal is to describe them as members of certain parametric families of complex Hadamard matrices of order $8$. It turns out that all matrices found in Proposition \ref{myCount} appear as members of $3$ distinct well-known parametric families of complex Hadamard matrices that we will now present. For each of these families we will give the classes of quaternary Hadamard matrices among them.

The first infinite family we have here is the $5$-parameter generalized Fourier family, stemming from the $8\times 8$ Fourier matrix as the following:
\[F_8^{(5)}(a,b,c,d,e)=\left[
\begin{array}{cccccccc}
 1 & 1 & 1 & 1 & 1 & 1 & 1 & 1 \\
 1 & a & b & c & -1 & -a & -b & -c \\
 1 & d & -1 & -d & 1 & d & -1 & -d \\
 1 & e & -b & -\overline{a}c e & -1 & -e & b & \overline{a}c e \\
 1 & -1 & 1 & -1 & 1 & -1 & 1 & -1 \\
 1 & -a & b & -c & -1 & a & -b & c \\
 1 & -d & -1 & d & 1 & -d & -1 & d \\
 1 & -e & -b & \overline{a}c e & -1 & e & b & -\overline{a}c e
\end{array}
\right].\]
We have chosen the nonstandard parametrization here to emphasize the fact that one can obtain $BH(4,8)$ matrices stemming from the matrix above if and only if the values of $a,b,c,d$ and $e$ are fourth roots of unity (see Remark \ref{myRem}).
\begin{rem}
As $F_8^{(5)}(a,b,c,d,e)$ features two rows and columns with entries $\pm1$ it is an infinite family of \emph{jacket matrices} (cf.\ \cite{lv}). Note that $F_8^{(5)}(a w,b w^2,c w^3,d w^2,e w^3)$ where $w=\mathbf{e}^{\frac{2\pi\mathbf{i}}{8}}$ is an eighth root of unity coincides with the matrix listed in \cite{karol}, and hence is of Di\c{t}\u{a}-type. The matrix corresponding to $F_8^{(5)}(-1/t^2,\mathbf{i},-\mathbf{i}/t^2,-\mathbf{i},-\mathbf{i}/t^2)$ is equivalent to an infinite family of circulant matrices containing in particular Horadam's matrix $K_4(\mathbf{i})$ (cf.\ \cite[p.\ $88$]{horadam}).
\end{rem}
After evaluating the matrix $F_8^{(5)}(a,b,c,d,e)$ at the $4^5=1024$ possible quintuples, we obtained the following
\begin{proposition}\label{myProp}
The family $F_8^{(5)}(a,b,c,d,e)$ contains $8$ inequivalent $BH(4,8)$ matrices: the symmetric matrices $F_8^{(5)}(1,1,1,1,1)$, $F_8^{(5)}(1,\mathbf{i},\mathbf{i},\mathbf{i},\mathbf{i})$, $F_8^{(5)}(\mathbf{i},1,\mathbf{i},1,\mathbf{i})$, $F_8^{(5)}(1,1,\mathbf{i},1,\mathbf{i})$ and further the matrices $F_8^{(5)}(1,1,1,1,\mathbf{i})$, $F_8^{(5)}(1,1,\mathbf{i},\mathbf{i},\mathbf{i})$, and their transpose.
\end{proposition}
To prove that the matrices in Proposition \ref{myProp} are indeed inequivalent, we need to recall a more powerful invariant than the Haagerup set. The \emph{fingerprint} of a complex Hadamard matrix $H$ of order $n\geq 4$ was introduced in \cite{finger} and reads
\[\Phi(H):=\left\{\left\{(v_i(d),m_i(d)):i\in I(d)\right\} : d=2,\hdots,\left\lfloor n/2\right\rfloor\right\},\]
where $I(d)$ is a finite index set counting the number of different moduli $v_i(d)$ taken by the $d\times d$ minors of $H$, the multiplicity of these values being denoted by $m_i(d)$, respectively. Note that $\Phi(H)$ is an ordered set. For instance, the fingerprint of the $8\times 8$ real Hadamard matrix is given by $\Phi(F_8^{(5)}(1,1,1,1,1))=$
\[\left\{\left\{(0, 336), (2, 448)\right\}, \left\{(0, 1344), (4, 1792)\right\},\left\{(0, 1428), (8,3136), (16, 336)\right\}\right\}.\]
That is, for $d=2$, $336$ minors are $0$ and $448$ have modulus $2$, \emph{etc}. In theory, the fingerprint is a more sophisticated invariant than the Haagerup set, because instead of looking at the $2\times 2$ minors only it takes into account the higher order minors as well, therefore capturing some of the global properties of complex Hadamard matrices. One of the limitations of the fingerprint is, however, that it cannot distinguish a matrix from its transpose.

Now with the aid of this invariant one readily verifies that the $15$ matrices appearing in Proposition \ref{myCount} correspond to $10$ different fingerprints, that is each matrix can be identified, up to transposition, with its own. It is straightforward to check that the $8$ matrices appearing in Proposition \ref{myProp} have $6$ different fingerprints, and as the family $F_8^{(5)}(a,b,c,d,e)$ is \emph{self-cognate} (cf.\ \cite{karol}) it follows that the transpose of the non-symmetric matrices are members of the family as well. Indeed, we have $\left(F_8^{(5)}(1,1,1,1,\mathbf{i})\right)^T=F_8^{(5)}(1,1,\mathbf{i},1,1)$ and $\left(F_8^{(5)}(1,1,\mathbf{i},\mathbf{i},\mathbf{i})\right)^T=F_8^{(5)}(1,\mathbf{i},\mathbf{i},1,\mathbf{i})$.

The list of the $BH(4,8)$ matrices is not complete yet, as in \cite{MRS} another family of $8\times 8$ complex Hadamard matrices were obtained from tiling abelian groups:
\[S_8^{(4)}(a,b,c,d)=\left[
\begin{array}{cccccccc}
 1 & 1 & 1 & 1 & 1 & 1 & 1 & 1 \\
 1 & d & -d & -d & -1 &  c d & - c d & d \\
 1 &  a\overline{d} &  b\overline{d} & - b\overline{d} & 1 & -1 & -1 & - a\overline{d} \\
 1 &  a & - b &  b & -1 & - c d &  c d & - a \\
 1 & -1 & - b\overline{d} &  b\overline{d} & 1 &  c & - c & -1 \\
 1 & -d &  b & - b & -1 & d & d & -d \\
 1 & - a\overline{d} & -1 & -1 & 1 & - c &  c &  a\overline{d} \\
 1 & - a & d & d & -1 & -d & -d &  a
\end{array}
\right].\]
It was shown in \cite{MRS} that the matrix $S_8^{(4)}(\mathbf{i},\mathbf{i},\mathbf{i},1)$ is not of Di\c{t}\u{a}-type, and hence a small neighborhood around it avoids the family $F_8^{(5)}(a,b,c,d,e)$ completely. The matrix corresponding to $S_8^{(4)}(1,1,\mathbf{i},\mathbf{i})$ is equivalent to the ``jacket conference matrix'' $J_8$ \cite{lv}. By evaluating the matrix $S_8^{(4)}(a,b,c,d)$ at the fourth roots of unity we found the following
\begin{proposition}\label{myProp2}
The family $S_8^{(4)}(a,b,c,d)$ and its transpose contain $8$ inequivalent $BH(4,8)$ matrices: the real Hadamard matrix $S_8^{(4)}(1,1,1,1)$ and the matrix $S_8^{(4)}(1,1,\mathbf{i},\mathbf{i})$, which are equivalent to a symmetric matrix, and further the matrices $S_8^{(4)}(1,1,1,\mathbf{i})$, $S_8^{(4)}(1,\mathbf{i},1,\mathbf{i})$, $S_8^{(4)}(1,1,\mathbf{i},1)$, and their transpose.
\end{proposition}
By analyzing the fingerprint of the obtained matrices, one can see that the following equivalences hold: $S_8^{(4)}(1,1,1,1)\sim F_8^{(5)}(1,1,1,1,1)$ and $\left(S_8^{(4)}(1,1,\mathbf{i},1)\right)^T \sim F_8^{(5)}(1,1,1,1,\mathbf{i})$, and therefore a further family is required to describe all of the $15$ $BH(4,8)$ matrices.

The following matrix was constructed from MUBs of order $4$ by Di\c{t}\u{a} \cite{dita3} very recently:
\[D_{8B}^{(5)}(a,b,c,d,e)=\left[
\begin{array}{cccccccc}
 1 & 1 & 1 & 1 & 1 & 1 & 1 & 1 \\
 1 & a & -a & d & -d & -a & a & -1 \\
 1 & b & b\overline{c}e & -d & d & -b\overline{c}e & -b & -1 \\
 1 & c & -e & -1 & -1 & e & -c & 1 \\
 1 & -c & e & -1 & -1 & -e & c & 1 \\
 1 & -b & -b\overline{c}e & -d & d & b\overline{c}e & b & -1 \\
 1 & -a & a & d & -d & a & -a & -1 \\
 1 & -1 & -1 & 1 & 1 & -1 & -1 & 1
\end{array}
\right].\]
We will see shortly that the family $D_{8B}^{(5)}(a,b,c,d,e)$ is essentially different from the families $F_8^{(5)}(a,b,c,d,e)$ and $S_8^{(4)}(a,b,c,d)$, as it accounts for most of the $BH(4,8)$ matrices. We remark here that the family $\left(D_{8B}^{(5)}(a,b,c,1,c)\right)^T$ is equivalent to the $3$-parameter matrix $P_8$ reported in \cite{IEEE}. We have the following
\begin{proposition}\label{myProp3}
The family $D_{8B}^{(5)}(a,b,c,d,e)$ together with its transpose contain $11$ inequivalent $BH(4,8)$ matrices, namely the real Hadamard matrix $D_{8B}^{(5)}(1,1,1,1,1)$, $D_{8B}^{(5)}(1,1,\mathbf{i},1,1)$ and the matrix $D_{8B}^{(5)}(1,\mathbf{i},\mathbf{i},1,\mathbf{i})$, which are equivalent to symmetric matrices, and the matrices $D_{8B}^{(5)}(1,1,1,1,\mathbf{i})$, $D_{8B}^{(5)}(1,1,1,\mathbf{i},\mathbf{i})$, $D_{8B}^{(5)}(1,1,\mathbf{i},\mathbf{i},1)$, $D_{8B}^{(5)}(1,\mathbf{i},\mathbf{i},\mathbf{i},\mathbf{i})$, and their transpose.
\end{proposition}
In particular, the matrix $D_{8B}^{(5)}(1,1,\mathbf{i},\mathbf{i},1)$ is not a member of any of the families $F_8^{(5)}(a,b,c,d,e)$, $S_8^{(4)}(a,b,c,d)$ or $\left(S_8^{(4)}(a,b,c,d)\right)^T$. We are ready to state the main result concerning $BH(4,8)$ matrices.
\begin{theorem}
There are precisely $15$ inequivalent $BH(4,8)$ matrices, which can be described through the parametric families $F_8^{(5)}(a,b,c,d,e), S_8^{(4)}(a,b,c,d)$ and $D_{8B}^{(5)}(a,b,c,d,e)$. Five of these matrices are equivalent to a symmetric matrix, namely $F_8^{(5)}(1,1,1,1,1)$, $F_8^{(5)}(1,\mathbf{i},\mathbf{i},\mathbf{i},\mathbf{i})$, $F_8^{(5)}(\mathbf{i},1,\mathbf{i},1,\mathbf{i})$, $F_8^{(5)}(1,1,\mathbf{i},1,\mathbf{i})$, and $S_8^{(4)}(1,1,\mathbf{i},\mathbf{i})$. The other $10$ matrices are given by $F_8^{(5)}(1,1,1,1,\mathbf{i})$, $F_8^{(5)}(1,1,\mathbf{i},\mathbf{i},\mathbf{i})$, $S_8^{(4)}(1,1,1,\mathbf{i})$, $S_8^{(4)}(1,\mathbf{i},1,\mathbf{i})$, $D_{8B}^{(5)}(1,1,\mathbf{i},\mathbf{i},1)$, and their transpose. These matrices, up to transposition, can be distinguished by their fingerprint.
\end{theorem}
Finally, we would like to point out that another well-known invariant of $BH(4,n)$ matrices, the \emph{Smith normal form} (SNF, cf.\ \cite{smf}) can detect $6$ equivalence classes only. In particular, the inequivalent matrices $F_8^{(5)}(1,\mathbf{i},\mathbf{i},\mathbf{i},\mathbf{i})$ and $S_8^{(4)}(1,1,\mathbf{i},\mathbf{i})$ share a common SNF. This is not surprising, as this invariant does not take into account the \emph{distribution} of the minors of complex Hadamard matrices.
\section{concluding remarks}
In this paper we have described all quaternary complex Hadamard matrices up to order $8$. The number of equivalence classes for $n=1,2,4,6,8$ read $1,1,2,1,15$, respectively. By utilizing continuous parametric families we could present the inequivalent matrices in a compact and elegant form, which allows one to study, display, and store these matrices efficiently. We remark here that the classification of the $n\times n$ \emph{generalized} Hadamard matrices over groups of order $4$ has been completed up to $n=16$ very recently \cite{hlt}, and we intend to continue this work and attempt to describe in a similar fashion the remaining $BH(4,n)$ matrices of size $n=10,12$ or even higher orders too in a forthcoming paper.

\section*{Acknowledgements} The author would like to thank the referees very much for their valuable comments and suggestions.

\end{document}